\documentclass{article}

\usepackage{amsmath}
\usepackage{amsfonts}
\title{\LARGE \textbf{On Interval Non-Edge-Colorable Eulerian Multigraphs}}
\author{Petros A. Petrosyan\\ \\Department of Informatics and Applied Mathematics,\\
Yerevan State University, 0025, Armenia\\
Institute for Informatics and Automation Problems,\\
National Academy of Sciences, 0014, Armenia\\ E-mail:
pet\_petros@ipia.sci.am}

\begin{document}
\textheight = 18.1cm \textwidth = 11.4cm \maketitle

\begin{abstract}
An edge-coloring of a multigraph $G$ with colors $1,\ldots,t$ is
called an interval $t$-coloring if all colors are used, and the
colors of edges incident to any vertex of $G$ are distinct and form
an interval of integers. In this note, we show that all Eulerian
multigraphs with an odd number of edges have no interval coloring.
We also give some methods for constructing of interval
non-edge-colorable Eulerian multigraphs.\\
\end{abstract}

\section{Introduction}

In this note we consider graphs which are finite, undirected, and
have no loops or multiple edges and multigraphs which may contain
multiple edges but no loops. Let $V(G)$ and $E(G)$ denote the sets
of vertices and edges of a multigraph $G$, respectively. The degree
of a vertex $v\in V(G)$ is denoted by $d_{G}(v)$, the maximum degree
of $G$ by $\Delta (G)$, and the chromatic index of $G$ by
$\chi^{\prime }\left(G\right)$. A multigraph $G$ is Eulerian if it
has a closed trail containing every edge of $G$. The terms and
concepts that we do not define can be found in \cite{West}.

A proper edge-coloring of a multigraph $G$ is a coloring of the
edges of $G$ such that no two adjacent edges receive the same color.
If $\alpha $ is a proper edge-coloring of $G$ and $v\in V(G)$, then
$S\left(v,\alpha \right)$ denotes the set of colors of edges
incident to $v$. A proper edge-coloring of a multigraph $G$ with
colors $1,\ldots,t$ is called an interval $t$-coloring if all colors
are used, and for any vertex $v$ of $G$, the set $S\left(v,\alpha
\right)$ is an interval of integers. A multigraph $G$ is interval
colorable if it has an interval $t$-coloring for some positive
integer $t$. The set of all interval colorable multigraphs is
denoted by $\mathfrak{N}$.

The concept of interval edge-coloring of multigraphs was introduced
by Asratian and Kamalian \cite{AsrKam1}. In \cite{AsrKam1,AsrKam2},
they proved the following result.\\

\noindent\textbf{Theorem 1.} If $G$ is a multigraph and $G\in
\mathfrak{N}$, then $\chi^{\prime }\left(G\right)=\Delta (G)$.
Moreover, if $G$ is a regular multigraph, then $G\in \mathfrak{N}$
if and only if $\chi^{\prime }\left(G\right)=\Delta (G)$.\\

Some results on interval edge-colorings of multigraphs were obtained
in \cite{Kam}. In \cite{PetKhach}, the authors described some
methods for constructing of interval non-edge-colorable bipartite
graphs and multigraphs.

In this note we show that all Eulerian multigraphs with an odd
number of edges have no interval coloring. We also give some methods
for constructing of interval non-edge-colorable Eulerian multigraphs.\\

\section{Results}

Let $G$ be a multigraph. For any $e\in E(G)$, by $G_{e}$ we denote
the multigraph obtained from $G$ by subdividing the edge $e$. For a
multigraph $G$, we define a multigraph $G^{\star}$ as follows:
\begin{center}
$V(G^{\star})= V(G)\cup \{u\}$, $u\notin V(G)$,\\
$V(G^{\star})= E(G)\cup \{uv:v\in V(G)~and~d_{G}(v)~is~odd\}$.
\end{center}
For a graph $G$, by $L(G)$ we denote the line graph of the graph
$G$.\\

We also need a classical result on Eulerian multigraphs.\\

\noindent\textbf{Euler's Theorem.} (\cite{Euler}) A connected
multigraph $G$ is Eulerian if and only if every vertex of $G$ has an even degree.\\

Now we can prove our result.\\

\noindent\textbf{Theorem 2.} If $G$ is an Eulerian multigraph and
$\vert E(G)\vert$ is odd, then $G\notin
\mathfrak{N}$.\\

\noindent\textbf{Proof} Suppose, to the contrary, that $G$ has an
interval $t$-coloring $\alpha$ for some $t$. Since $G$ is an
Eulerian multigraph, $G$ is connected and $d_{G}(v)$ is even for any
$v\in V(G)$, by Euler's Theorem. Since $\alpha$ is an interval
coloring and all degrees of vertices of $G$ are even, we have that
for any $v\in V(G)$, the set $S\left(v,\alpha\right)$ contains
exactly $\frac{d_{G}(v)}{2}$ even colors and $\frac{d_{G}(v)}{2}$
odd colors. Now let $m_{odd}$ be the number of odd colors in the
coloring $\alpha$. By Handshaking lemma, we obtain
$m_{odd}=\frac{1}{2}\sum\limits_{v\in
V(G)}\frac{d_{G}(v)}{2}=\frac{\vert E(G)\vert}{2}$. Thus $\vert
E(G)\vert$ is even, which is a contradiction. $\square$\\

\noindent\textbf{Corollary 1.} If $G$ is an Eulerian multigraph and
$G\in \mathfrak{N}$, then $\vert E(G)\vert$ is even.\\

Let us note that there are Eulerian graphs with an even number of
edges that have no interval coloring. For example, the complete
graph $K_{5}$ has no interval coloring. On the other hand, there are
many Eulerian graphs with an even number of edges that have an
interval coloring. In
\cite{Jaeger}, Jaeger proved the following result.\\

\noindent\textbf{Theorem 3.} If $G$ is a connected $r$-regular graph
($r\geq 2$), $\chi^{\prime }\left(G\right)=r$ and $\vert E(G)\vert$
is even, then
$\chi^{\prime }\left(L(G)\right)=2r-2$.\\

Since $G$ is a connected $r$-regular graph ($r\geq 2$) and $\vert
E(G)\vert$ is even, we have that $L(G)$ is a connected
$(2r-2)$-regular graph with an even number of edges. Moreover, by
Theorems 1 and 3 and Euler's Theorem,
 we obtain the following\\

\noindent\textbf{Corollary 2.} If $G$ is a connected $r$-regular
($r\geq 2$) graph with an even number of edges and $G\in
\mathfrak{N}$, then
$L(G)$ is an Eulerian graph with an even number of edges and $L(G)\in \mathfrak{N}$.\\

Let us note that Theorem 2 also gives some methods for constructing
of interval non-edge-colorable Eulerian multigraphs from interval
colorable multigraphs.\\

\noindent\textbf{Corollary 3.} If $G$ is an Eulerian multigraph and
$G\in \mathfrak{N}$, then for each $e\in E(G)$, $G_{e}\notin
\mathfrak{N}$.\\

\noindent\textbf{Corollary 4.} If $G$ is a connected multigraph with
an odd number of edges and $G\in \mathfrak{N}$, then
$G^{\star}\notin \mathfrak{N}$.\\

\end{document}